\documentclass[a4paper,11pt]{article}

\usepackage{amsmath,amssymb,bm,ulem,amsthm}
\usepackage[dvipdfmx]{graphicx}
\usepackage{here}
\usepackage{url}
\usepackage{comment}

\usepackage[dvipdfmx]{graphicx}

\theoremstyle{plain}
\newtheorem{lemma}{Lemma}
\newtheorem{prop}[lemma]{Proposition}
\newtheorem{thm}[lemma]{Theorem}

\title{Computation of Vector-Valued Invariants for a Finite Complex Reflection Group}
\author{ A.K.M. Selim Reza\footnote{Corresponding author},  Manabu Oura, Masashi Kosuda and Shoyu Nagaoka }

\begin{document}
\date{\empty} 
\maketitle

\textbf{Abstract}. We consider the complex reflection group \( \mathcal{G}  \), identified as No. 8 in the Shephard-Todd classification. In this paper, we present computations of the vector-valued invariants associated with various representations of \( \mathcal{G}  \). Additionally, we investigate the structure of the corresponding invariant rings.\\



\section{Introduction}
\label{Intro}
In this study, we present computations of the vector-valued invariants of \( \mathcal{G}  \) with respect to a representation that will be defined later. We begin by reviewing the invariant theory of a finite group.

Let \( G \) be a finite subgroup of \( GL(n, \mathbf{C}) \). We denote a set of variables by \( \bm{x} = \begin{pmatrix} x_1 \\ x_2 \\ \vdots \\ x_n \end{pmatrix} \). The group \( G \) acts on \( \mathbf{C}[\bm{x}] = \mathbf{C}[x_1, x_2, \dots, x_n] \) as follows:
\[
g . f(\bm{x}) = f(g \bm{x}), \quad g \in G
\]
where \( g \bm{x} \) represents the usual product of an \( n \times n \) matrix \( g \) and the column vector \( \bm{x} \). A function \( f(\bm{x}) \) is said to be an invariant of \( G \) if \( g . f(\bm{x}) = f(\bm{x}) \) for all \( g \in G \). The set of all invariants of \( G \) is denoted by \( \mathbf{C}[\bm{x}]^G \), which is known as the invariant ring of \( G \).

For the purpose of this paper, we focus on the specific group \( \mathcal{G}  \). We briefly recall the properties of \( \mathcal{G}  \) and its definition. Let
\[
\mathcal{G}  = \langle T, D \rangle,
\]
where
\[
T = \frac{1+i}{2} \begin{pmatrix} 1 & 1 \\ 1 & -1 \end{pmatrix}, \quad D = \begin{pmatrix} 1 & 0 \\ 0 & i \end{pmatrix}.
\]
This group has order 96 and belongs to the 37 classes of complex reflection groups. Specifically, \( \mathcal{G}  \) is the 8th group in the list of Shephard-Todd classification \cite{shephard-todd}. Since \( \mathcal{G}  \) consists of \( 2 \times 2 \) matrices, we adopt the notation \( \bm{x} = \begin{pmatrix} x \\ y \end{pmatrix} \) instead of \( \bm{x} = \begin{pmatrix} x_1 \\ x_2 \end{pmatrix} \).

It is known that the invariant ring \( \mathbf{C}[\bm{x}]^{\mathcal{G} } \) is generated by two algebraically independent polynomials:
\[
\theta = x^8 + 14x^4y^4 + y^8, \quad \varphi = x^{12} - 33x^8y^4 - 33x^4y^8 + y^{12}.
\]
The dimension formula for this ring is given by
\[
\frac{1}{(1 - t^8)(1 - t^{12})}.
\]

The group \( \mathcal{G}  \) and its invariant ring are of particular interest in number theory. To elaborate on this, for \( a, b \in \{0,1\} \), we define the theta functions
\[
\theta_{ab}(\tau) = \sum_{n \in \mathbf{Z}} \exp \left( 2\pi \sqrt{-1} \left( \frac{\tau}{2} \left( n + \frac{a}{2} \right)^2 + \left( n + \frac{a}{2} \right) \frac{b}{2} \right) \right),
\]
where \( \tau = x + \sqrt{-1} y \), with \( x \in \mathbf{R} \) and \( y \in \mathbf{R}_{> 0} \). The map
\[
x \mapsto \theta_{00}(2\tau), \quad y \mapsto \theta_{10}(2\tau)
\]
induces an isomorphism between the invariant ring \( \mathbf{C}[\bm{x}]^{\mathcal{G} } \) and the ring of modular forms for \( SL(2, \mathbf{Z}) \), as shown in \cite{broueenguehard, rungecodes}.\\

In the earlier work \cite{kosuda-oura}, second and third named authors investigated the structure of the centralizer ring of the tensor products of \( \mathcal{G}  \). In the present paper, we explore our study by examining further properties of \( \mathcal{G}  \) from the perspective of invariant theory. All computations in this work are performed using SageMath \cite{sagemath} and Maple.

\section{Preliminaries}
Let \( G \) be a finite subgroup of \( GL(n, \mathbf{C}) \). We fix a representation of \( G \) of degree \( m \), \( \rho: G \to GL(m, \mathbf{C}) \). Corresponding to the degree \( m \) of the fixed representation, we define
\[
F(\bm{x}) = \begin{pmatrix} f_1(\bm{x}) \\ f_2(\bm{x}) \\ \vdots \\ f_m(\bm{x}) \end{pmatrix},
\]
where \( f_i(\bm{x}) \)'s are homogeneous polynomials of the same degree in \( \mathbf{C}[\bm{x}] \). For \( g \in G \), we define the action of \( G \) on \( F(\bm{x}) \) as
\[
g . F(\bm{x}) = \begin{pmatrix} g . f_1(\bm{x}) \\ g . f_2(\bm{x}) \\ \vdots \\ g . f_m(\bm{x}) \end{pmatrix} = \begin{pmatrix} f_1(g \bm{x}) \\ f_2(g \bm{x}) \\ \vdots \\ f_m(g \bm{x}) \end{pmatrix}.
\]

We observe that \( (gh) . F(\bm{x}) = F\left( (gh)\bm{x} \right) = h . (g . F(\bm{x})) \). We say that \( F(\bm{x}) \) is a vector-valued invariant of \( G \) with respect to \( \rho \) if

\begin{equation}
\label{vinv}
g \cdot F(\bm{x}) = \rho(g) F(\bm{x})
\end{equation}

for every \( g \in G \) (\textit{cf.} \cite{freitag}). Notice that \( \rho(g) F(\bm{x}) \) is the usual matrix product of an \( m \times m \) matrix \( \rho(g) \) and a column vector \( F(\bm{x}) \) of size \( m \).\\

The set \( M(\rho) \) of all vector-valued invariants of \( G \) with respect to \( \rho \) forms a \( \mathbf{C}[\bm{x}]^G \)-module. The set of elements of degree \( k \) in \( M(\rho) \) forms a finite-dimensional vector space over \( \mathbf{C} \).


\begin{prop}\label{reynolds}
For any $F(\bm{x})$, 
\[
\frac{1}{|G|} \sum_{g\in G} \rho(g)^{-1} F(g\bm{x})
\]
is a vector valued invariant of $G$ with respect to $\rho$.
\end{prop}

\textit{Proof.} 
In order to show that $F'(\bm{x})=\frac{1}{|G|} \sum_{g\in G} \rho(g)^{-1} F(g\bm{x})$ is a vector valued invariant with respect to a representation $\rho$,
we need to prove
\[
\rho(h)^{-1}F'(h\bm{x})=F'(\bm{x})
\]
for every $h\in G$.
Take $h\in G$. We have
\begin{align*}
\rho(h)^{-1}\left( \frac{1}{|G|} \sum_{g\in G} \rho(g)^{-1}F(gh\bm{x})\right) &=
\frac{1}{|G|} \sum_{g\in G} \rho(h^{-1}g^{-1})F(gh\bm{x})\\
&=\frac{1}{|G|} \sum_{g\in G} \rho((gh)^{-1})F(gh\bm{x})\\
&=\frac{1}{|G|} \sum_{g\in G} \rho(gh)^{-1}F(gh\bm{x})\\
&=\frac{1}{|G|} \sum_{g\in G} \rho(g)^{-1}F(g\bm{x}).
\end{align*}
Thus, we have shown that \( F'(\bm{x}) \) is a vector-valued invariant with respect to \( \rho \), proving Proposition \ref{reynolds}.

Let $\mathcal{G}  \subset \mathrm{GL}_2(\mathbb{C})$ be the group defined in Section~\ref{Intro}. For any finite-dimensional representation $\rho : \mathcal{G}  \to \mathrm{GL}_m(\mathbb{C})$, we define 
\[
M(\rho) \;=\; \big\{\, F \in \mathbb{C}[x,y]^m \ \big|\ F(g\cdot (x,y)) \;=\; \rho(g)\,F(x,y) \ \text{for all } g \in \mathcal{G}  \,\big\},
\]
the $\mathbb{C}[x,y]^{\mathcal{G} }$–module of vector–valued polynomials whose components transform according to the representation~$\rho$. The dimension of $M(\rho)$ can be obtained using the equivariant Molien formula (cf.~\cite{gatermann-1996}, \cite{worfolk-1994}).
\begin{equation}\label{eq:equivariant-molien}
H_{M(\rho)}(t) \;=\; \frac{1}{|\mathcal{G} |} \sum_{g\in \mathcal{G} } 
\frac{\mathrm{tr}\,\rho(g^{-1})}{\det(I - t\,g)},
\end{equation}
where the determinant in the denominator is computed in the natural $2\times 2$ representation of $\mathcal{G} $ on $(x,y)$, while the trace in the numerator is taken in the chosen representation~$\rho$.

In this case, $\mathbb{C}[x,y]^{\mathcal{G} } \cong \mathbb{C}[\theta,\phi]$ with $\deg \theta = 8$ and $\deg \phi = 12$, so the denominator of $H_{M(\rho)}(t)$ always factors as
\[
(1 - t^{8})(1 - t^{12}).
\]

\medskip
Kosuda \& Oura~\cite{kosuda-oura} explicitly compute all irreducible 
representations $\rho_i$ (for $i = 1, 2, \dots, 16$) of the group $\mathcal{G}$, giving each representation by explicit matrix generators in terms of the chosen presentation of $\mathcal{G}$. In particular, their results show 
that, the irreducible representations satisfy
\[
\begin{array}{lll}
\rho_2 = \rho_3 \otimes \rho_3, & 
\rho_4 = \rho_2 \otimes \rho_3, & 
\rho_6 = \rho_3 \otimes \rho_5, \\[2pt]
\rho_7 = \rho_3 \otimes \rho_{10}, & 
\rho_8 = \rho_2 \otimes \rho_{10}, & 
\rho_9 = \rho_4 \otimes \rho_{10}, \\[2pt]
\rho_{11} = \rho_3 \otimes \rho_{13}, & 
\rho_{12} = \rho_2 \otimes \rho_{13}, & 
\rho_{14} = \rho_4 \otimes \rho_{13}, \\[2pt]
\rho_{16} = \rho_3 \otimes \rho_{15}. & & 
\end{array}
\]
together with the fundamental representations 
\[
\rho_1, \quad \rho_3, \quad \rho_5, \quad \rho_{10}, \quad \rho_{13}, \quad \rho_{15}.
\]
from which all others can be obtained by tensor products.\\

Let
\[
G_i = \langle \rho_i(T), \rho_i(D) \rangle.
\]
The order of each \( G_i \) is given in \cite{kosuda-oura}. For the convenience of the readers, we reproduce the results here.

\begin{prop}
The order of \( G_i \) for \( i = 1, 3, 5, 10, 13, 15 \) is
\[
1, \ 4, \ 6, \ 96, \ 48, \ 96
\]
respectively.
\end{prop}

In next section we use these six representations to investigate the structure of the vector-valued invariants as well as their's corresponding dimension formula via the equivariant Molien formula. 

\bigskip

\section{Results}

We use the set \( \bm{x} = \begin{pmatrix} x \\ y \end{pmatrix} \) of variables instead of \( \bm{x} = \begin{pmatrix} x_1 \\ x_2 \end{pmatrix} \). In this case, we have \( \mathbf{C}[\bm{x}] = \mathbf{C}[x, y] \).\\

\noindent\underline{The case \( \rho_1 \)}

\noindent This is a one-dimensional (trivial) representation which maps each element of the group to 1. We denote that of \( H_1 \) by \( (\rho_1, V_1) \). Thus, we have:
\[
\rho_1(T) = 1, \quad \rho_1(D) = 1.
\]
This is a standard invariant, as for every \( g = \begin{pmatrix} a & b \\ c & d \end{pmatrix} \in \mathcal{G}  \),
\[
g . f(\bm{x}) = f(\bm{x}).
\]
Therefore, it holds that \( M(\rho_1) = \mathbf{C}[\bm{x}]^{\mathcal{G} } \). The invariant ring \( \mathbf{C}[\bm{x}]^{\mathcal{G} } \) is generated by two algebraically independent polynomials, \( \theta \) and \( \varphi \). The dimension formula for this ring is
\[
\frac{1}{(1 - t^8)(1 - t^{12})} = 1 + t^8 + t^{12} + t^{16} + \dots.
\]
This case is well-known in the literature.

\bigskip
\noindent\underline{The case \( \rho_3 \)}

\noindent The determinant which maps $T$ and $D$ to $-i$ and $i$ respectively also gives another one-dimensional irreducible representation. We call it $(\rho_3,V_3)$.
Thus, we have:
\[
\rho_3(T) = -i = \det(T), \quad \rho_3(D) = i =\det( D)
\]

\noindent We observe that the extended group \(G_1= \langle \mathcal{G} , \eta_8 \rangle \) has order 192 and is a complex reflection group of No.9 in the Shephard-Todd classification. The invariance property of \( G_1 \) with respect to the determinant representation is discussed in \cite{bachoc}.

\noindent Using equation~(\ref{eq:equivariant-molien}), we get the dimension formula for \( M(\rho_3) \) as
\begin{align*}
    H_{M(\rho_{3})}(t) 
    &= \frac{t^6}{(1 - t^8)(1 - t^{12})} \\
    &= t^6 + t^{14} + t^{18} + t^{22} + t^{26} + \dots.
\end{align*}

\medskip
\noindent\underline{The case \( \rho_5 \)}

\noindent The representation \( \rho_5 \) is a two-dimensional irreducible representation where the group elements \( T \) and \( D \) act on the vector space say, \( V_5 \) via the matrices:
    \[
    \rho_5(T) = \frac{-1}{2} \begin{pmatrix} 1 & 1 \\ 3 & -1 \end{pmatrix}, \quad \rho_5(D) = \text{diag}(1, -1).
    \]

\noindent By equation~(\ref{eq:equivariant-molien}), the dimension of \( M(\rho_5) \) is
\begin{align*}
    H_{M(\rho_{5})}(t) 
    &= \frac{t^4}{(1 - t^4)(1 - t^{12})} \\
    &= \frac{t^4 + t^8}{(1 - t^8)(1 - t^{12})}\\ 
    &= t^4 + t^8 + t^{12} + 2t^{16} + 2t^{20} + \dots
\end{align*}

\medskip
\noindent\underline{The case \( \rho_{10} \)}

\noindent The representation \( \rho_{10} \) is a two-dimensional representation that maps \( T \) and \( D \) to their defining matrices:
    \[
    \rho_{10}(T) = T, \quad \rho_{10}(D) = D,
    \]
where these matrices describe how \( T \) and \( D \) act on the vector space say, \( V_{10} \).

\noindent Applying equation~(\ref{eq:equivariant-molien}), we obtain the dimension formula for
\( M(\rho_5) \):
\begin{align*}
    H_{M(\rho_{10})}(t) 
    &= \frac{t}{(1 - t^4)(1 - t^{12})} \\
    &= \frac{t + t^5}{(1 - t^8)(1 - t^{12})}\\ 
    &=  t + t^5 + t^9 + 2t^{13} + 2t^{17} + 2t^{21} + \dots
\end{align*}

\medskip
\noindent\underline{The case \( \rho_{13} \)}

\noindent The representation \( \rho_{13} \) is a three-dimensional irreducible representation, where the group elements \( T \) and \( D \) act on the vector space say, \( V_{13} \) via the matrices:
    \[
    \rho_{13}(T) = \frac{i}{2} \begin{pmatrix} 1 & 2 & 1 \\ 1 & 0 & -1 \\ 1 & -2 & 1 \end{pmatrix}, \quad \rho_{13}(D) = \text{diag}(1, i, -1).
    \]
By equation~(\ref{eq:equivariant-molien}), the dimension of \( M(\rho_{13}) \) is
\begin{align*}
    H_{M(\rho_{13})}(t) 
    &= \frac{t^2}{(1 - t^4)(1 - t^{8})} \\
    &= \frac{t^2 + t^6+t^{10}}{(1 - t^8)(1 - t^{12})}\\ 
    &=  t^2 + t^6 + t^{10} + 2t^{14} + 3t^{18} + 3t^{22} + \dots
\end{align*}
\medskip
\noindent\underline{The case \( \rho_{15} \)}

\noindent The representation \( \rho_{15} \) is a four-dimensional irreducible representation, where the group elements \( T \) and \( D \) act on the vector space \( V_{15} \) via the matrices:
    \[
    \rho_{15}(T) = \frac{-1 + i}{4} \begin{pmatrix} 1 & 3 & 3 & 1 \\ 1 & 1 & -1 & -1 \\ 1 & -1 & -1 & 1 \\ 1 & -3 & 3 & -1 \end{pmatrix}, \quad \rho_{15}(D) = \text{diag}(1, i, -1, -i).
    \]
By equation~(\ref{eq:equivariant-molien}), the dimension of \( M(\rho_{15}) \) is
\begin{align*}
    H_{M(\rho_{15})}(t) 
    &= \frac{t^3+t^{11}}{(1 - t^4)(1 - t^{12})} \\
    &= \frac{t^3 + t^7+t^{11}+t^{15}}{(1 - t^8)(1 - t^{12})}\\ 
    &=  t^3 + t^7 + 2t^{11} + 3t^{15} + 3t^{19} + 4t^{23} + \dots
\end{align*}

\noindent To get the corresponding vector valued invariants first we prepare the function \( F(x) \) with coefficients as parameters. We then apply the condition defined in the equation (\ref{vinv}) for vector-valued invariants with respect to \( \rho_i \) to \( F(x) \). By solving these conditions, we obtain a free \( \mathbf{C}[\bm{x}]^{\mathcal{G} } \)-module that is in \( M(\rho_i) \).

\bigskip

Thus, we have obtained the following theorem.

\begin{thm}
\begin{enumerate}
\item
$M(\rho_{3})$ is a free $\mathbf{C}[\bm{x}]^{\mathcal{G} }$ module generated by 
\[
-x^5y+xy^5
\]
whose dimension formula is 
\[
\frac{t^6}{(1-t^8)(1-t^{12})}=t^6+t^{14}+t^{18}+t^{22}+t^{26}+\dots.
\]

\item
$M(\rho_5)$ is a free $\mathbf{C}[\bm{x}]^{\mathcal{G} }$ module generated by 
\[
\begin{pmatrix} x^4+y^4 \\ 6x^2y^2\end{pmatrix},\
\begin{pmatrix}
-x^8+10x^4y^4-y^8\\ 12x^6y^2+12x^2y^6
\end{pmatrix}
\]
whose dimension formula is 
\[
\frac{t^4+t^8}{(1-t^8)(1-t^{12})}=t^4+t^8+t^{12}+2t^{16}+2t^{20}+\dots .
\]

\item
$M(\rho_{10})$ is a free $\mathbf{C}[\bm{x}]^{\mathcal{G} }$ module generated by 
\[
\begin{pmatrix} x \\ y \end{pmatrix},\
\begin{pmatrix}
-x^5+5xy^4\\ 5x^4y-y^5
\end{pmatrix}
\]
whose dimension formula is 
\[
\frac{t+t^5}{(1-t^8)(1-t^{12})}=t+t^5+t^9+2t^{13}+2t^{17}+2t^{21}+\dots.
\]

\item
$M(\rho_{13})$ is a free $\mathbf{C}[\bm{x}]^{\mathcal{G} }$ module generated by 
\[
\begin{pmatrix} x^2\\ xy\\ y^2\end{pmatrix},\
\begin{pmatrix} -x^6+5x^2y^4\\
2x^5y+2xy^5\\
5x^4y^2-y^6
\end{pmatrix},\
\begin{pmatrix}
-x^2(x^4-5y^4)^2\\
xy(5x^4-y^4)(x^4-5y^4)\\
-y^2(5x^4-y^4)^2
\end{pmatrix}
\]
whose dimension formula is

\[
\frac{t^2+t^6+t^{10}}{(1-t^8)(1-t^{12})}=t^2+t^6+t^{10}+2t^{14}+3t^{18}+3t^{22}+\dots.
\]

\item
$M(\rho_{15})$ is a free $\mathbf{C}[\bm{x}]^{\mathcal{G} }$ module generated by 
\[
\begin{pmatrix}
x^3\\ x^2y\\ xy^2\\ y^3 \end{pmatrix},\
\begin{pmatrix}
-x^7+5x^3y^4\\
x^6y+3x^2y^5\\
3x^5y^2+xy^6\\
5x^4y^3-y^7\end{pmatrix},\
\begin{pmatrix}
-6x^7y^4+6x^3y^8\\ 
-x^{10}y+x^2y^9\\
x^9y^2-xy^{10}\\
6x^8y^3-6x^4y^7
\end{pmatrix},\
\begin{pmatrix}
-24x^{11}y^4+24x^7y^8\\
 x^{14}y-7x^{10}y^5+3x^6y^9+3x^2y^{13}\\
  3x^{13}y^2+3x^9y^6-7x^5y^{10}+xy^{14}\\ 24x^8y^7-24x^4y^{11}
\end{pmatrix}
\]
whose dimension formula is 

\[
\frac{t^3+t^7+t^{11}+t^{15}}{(1-t^8)(1-t^{12})}=t^3+t^7+2t^{11}+3t^{15}+3t^{19}+4t^{23}+\dots.
\]
\end{enumerate}
\end{thm}

\bigskip

\textbf{Acknowledgements}.
We thank the reviewer for his insightful comment. The first named author was supported in part by funds from  Ministry of Education, Culture, Sports, Science and Technology (MEXT), Japan and the second named author was supported by JSPS KAKENHI Grant Numbers 24K06827 and 24K06644. 



Department of Mathematics, Khulna University of Engineering \& Technology, Khulna-9203, Bangladesh and
Graduate School of Natural Science and Technology, Kanazawa University, Ishikawa 920-1192, Japan

\textit{Email address}: \texttt{selim\_1992@math.kuet.ac.bd}

\bigskip

Faculty of Mathematics and Physics,
Institute of Science and Engineering,
Kanazawa University,
Kakuma-machi,
Ishikawa 920-1192,
Japan

\textit{Email address}: \texttt{oura@se.kanazawa-u.ac.jp}

\bigskip

Faculty of Engineering, University of Yamanashi, 400-8511, Japan

\textit{Email address}: \texttt{mkosuda@yamanashi.ac.jp}

\bigskip

Emeritus Professor of Kindai University, Osaka, 545-0001, Japan

\textit{Email address}: \texttt{shoyu1122.sn@gmail.com}

\end{document}